\documentclass{article}
\usepackage{arxiv}
\usepackage{hyperref}       
\usepackage{url}            
\usepackage{booktabs}       
\usepackage{amsfonts}       
\usepackage{nicefrac}       
\usepackage{microtype}      
\usepackage{lipsum}		
\usepackage{graphicx}
\usepackage{natbib}
\usepackage{doi}
\usepackage{dcolumn}
\usepackage{bm}
\usepackage{verbatim}
\usepackage{amsmath}
\usepackage{amssymb}
\newtheorem{lemma}{Lemma}
\newtheorem{theorem}[lemma]{Theorem}
\newtheorem{proposition}[lemma]{Proposition}

\newtheorem{proof}[lemma]{Proof}
\title{Spectral properties of the Laplacian of temporal networks following a constant block Jacobi model}

\author{{Zhana Kuncheva} \\ 
Data Science and Engineering\\
Optima Partners\\
London, UK\\
zhana.kuncheva@optimapartners.co.uk\\
\And
{Ognyan Kounchev}\\
Department of Mathematics and Computer Science\\
Bulgarian Academy of Sciences\\
Sofia, Bulgaria\\
kounchev@math.bas.bg}

\date{\today}
\begin{document}
\maketitle
\begin{abstract}
We study the behavior of the eigenvectors associated with the smallest eigenvalues
of the Laplacian matrix of temporal networks. We consider the multilayer representation of temporal networks, i.e. a set of networks linked through ordinal interconnected layers. We
analyze the Laplacian matrix, known as supra-Laplacian, constructed through the
supra-adjacency matrix associated with the multilayer formulation of temporal
networks, using a constant block Jacobi model which has closed-form solution. To do this, we assume that the inter-layer weights are perturbations of the
Kronecker sum of the separate adjacency matrices forming the temporal network.
Thus we investigate the properties of the eigenvectors associated with the
smallest eigenvalues (close to zero) of the supra-Laplacian matrix. Using
arguments of perturbation theory, we show that these eigenvectors can be
approximated by linear combinations of the zero eigenvectors of the individual
time layers.
This finding is crucial in reconsidering and generalizing the role of the Fielder vector in
supra-Laplacian matrices.

\end{abstract}

\maketitle

\section{Introduction}
In recent years, one of the major lines of research in complex network
analysis is the topological changes that occur in a network over time. A
sequence of networks with such a time-varying nature can be formalized as a
\textit{temporal network} \cite{Holme2012}. The \textit{multilayer
formulation} of temporal networks \cite{Kivel} is one way to consider the
interconnected topological structure changing over time: \textit{ordinal}
interconnections between layers determine how a given node in one layer and
its given counterparts in the previous and next time point layers are linked
and influence each other. The network analysis community has strong traditions
in using the spectral properties \cite{Moreno2013, Sol2013} of multilayer
networks for various purposes such as centrality measures
\cite{DeDomenico2013} or investigating diffusion processes \cite{Sol2013}.

One challenge associated with understanding the spectral properties of the
temporal networks is the lack of available tools that respect the fundamental
distinction between within-layer and inter-layer edges
\cite{Kivel,Taylor2015,DeDomenico2013a} when studying the spectral properties
of the Laplacian matrix $\mathcal{L}$ of temporal networks, known as
supra-Laplacian. A number of investigations were undertaken to show that the
inter-layer couplings in multilayer networks distort those spectral properties
and to explain the effect of different inter-layer weights over the
eigenvalues of the supra-Laplacian \cite{Moreno2013, Sol2013}. Up to our
knowledge, there is no work related to the understanding of the information
carried by the eigenvectors corresponding to the smallest eignevalues of the supra-Laplacian.

The spectral analysis on a network is nowadays understood as studying the
spectral properties of the various Laplacian matrices defined on the network.
In particular, for the so-called normalized Laplacian the most interesting are
usually the smallest eigenvalues and their eigenvectors.

For a Laplacian matrix, the eigenvector corresponding to the smallest
eigenvalue, $\lambda_{1}=0$, is constant or weighted by the node degrees if
the Laplacian is normalized \cite{Chung1996}. The eigenvector corresponding
to the smallest non-zero eigenvalue, known as the algebraic connectivity, is
in practice used for partitioning purposes \cite{Luo2002,Luxburg} and is
known as the Fiedler vector. In this article, we consider
\emph{slowly-changing temporal networks} which means that the adjacency
matrices forming the different time layers change relatively slowly~\cite{ENRIGHT201888}. The main objective
of the present paper is to draw a maximal profit of this important property for
the majority of temporal networks. In particular, for every temporal network, for a sufficiently small interval, we have this effect.

Further, we add inter-layer weights to the temporal network which may be
considered as \textit{perturbations} of the Kronecker sum of the separate
adjacency matrices forming the different time layers, and we consider the
Laplacian of the resulting matrix which is usually called
supra-Laplacian \cite{Kivel}. This point of view on the temporal networks,
allows us to find an approximate closed form solution of the eigenvectors
corresponding to the smallest eigenvalues of the supra-Laplacian. In
particular, by applying arguments from perturbation theory, we are able to
show that the eigenvectors corresponding to the smallest eigenvalues (of the
supra-Laplacian) are well approximated by the space of the perturbed
eigenvectors corresponding to all zero eigenvalues of the Laplacian matrices
corresponding to the networks of the separate time layers.

The paper is organised as follows: in Sec.~\ref{sec2}, we present the
construction of the temporal network following a constant block Jacobi model. This model appears in a natural way as a first order approximation to the slowly-changing temporal network, and enjoys a closed-form solution of the eigenvectors of the supra-Laplacian matrix; in Sec.~\ref{sec3} we
investigate the spectral properties of the supra-Laplacian and obtain an eigenvector solution of the reduced system;
Sec.~\ref{sec5} is devoted to identifying the smallest eigenvectors, which are
obtained by perturbation of the zero eigenvectors of the separate time layers,
and discussing the influence of density and number of layers on these
eigenvectors; finally we state the conclusions.

\section{Temporal network following constant block Jacobi model: notations and
definitions\label{sec2}}
A \textit{temporal network} is a set of networks in which edges and nodes vary
in time. In this work, we make the assumption that each node $i$ is present in
all layers. We use the notation $G^{t}$ for a layer in an ordered sequence of
$T$ networks $\mathcal{T}=\left\{  G^{1},G^{2},...,G^{T}\right\}  $ with
$G^{t}=\left(  V,A^{t}\right)  $ where $t\in\left\{  1,2,...,T\right\}  $ and
the number of nodes is $N,$ i.e. $N=\left\vert V\right\vert .$ Here $A^{t}$ is
a binary undirected and connected adjacency matrix. In order to use the multilayer framework for
representing a \textit{temporal network}, we consider the \textit{diagonal
ordinal coupling} of layers \cite{Kivel,Bassett2011,Mucha}, to define a new
supra-network $\widetilde{\mathcal{T}}$ . We define the coupling edges: we denote by
$\omega_{i}^{t,p}\in\mathbb{R}$ the value of the inter-layer edge weight
between node $i$ in different time layers $t$ and $p$. Our main assumption is
that only neighbouring layers may be connected, i.e. $\omega_{i}^{t,p}=0$ for
all layers $G^{t}$ and $G^{p}$, with $p\neq t-1$ and $p\neq t+1$. No other
edges between $G^{t}$ and $G^{p}$ exist for indices $t\neq p.$

As a result, the multilayer framework of the temporal network is expressed in
an $NT$-node single adjacency matrix $\mathcal{A}$ of size $NT\times NT$ which
is simply the adjacency matrix of the network $\widetilde{\mathcal{T}}$ ,
referred to as \emph{supra-adjacency matrix}. Clearly, the diagonal blocks of
$\mathcal{A}$ are the adjacency matrices $A^{t}$, and the off-diagonal blocks
are the inter-layer weight matrices $W^{t,p}=diag(\omega_{1}^{t,p},\omega
_{2}^{t,p},...,\omega_{N}^{t,p})$ if $p=t-1$ or $p=t+1.$ 

The usual within-layer degree of node $i$ in layer $G^{t}$ is defined as
$d_{i}^{t}:={\sum_{j=1}^{N}}A_{ij}^{t}$ while the multilayer node degree of
node $i$ in layer $G^{t}$ is $\mathfrak{d}_{i}^{t}:=d_{i}^{t}+\omega
_{i}^{t,t-1}+\omega_{i}^{t,t+1}$. Define the degree matrix $\mathcal{D}$ as
$\mathcal{D}:=\text{diag}\left(  \mathfrak{d}_{1}^{1},\mathfrak{d}_{2}^{1}%
,...,\mathfrak{d}_{N}^{1},\mathfrak{d}_{1}^{2},...,\mathfrak{d}_{N}%
^{2},...,\mathfrak{d}_{N}^{T}\right) $. The \textit{normalized supra-Laplacian}
$\mathcal{L}$ is defined as $\mathcal{L}:\mathfrak{=}\mathcal{D}^{-\frac{1}{2}}\left(\mathcal{D-A}\right)\mathcal{D}^{-\frac{1}{2}}$ \cite{Chung1996}.

The supra-adjacency matrix $\mathcal{A}^{0}$ with $0$ inter-layer weights and
its corresponding Laplacian matrix $\mathcal{L}^{0}$ are directly expressed as
a Kronecker sum:
\begin{equation}
\mathcal{A}^{0}:=\oplus_{t=1}^{T}A^{t}\longrightarrow\mathcal{L}^{0}%
=\oplus_{t=1}^{T}L^{t} \label{eq:tensorA}%
\end{equation}
where $L^{t}$ is the normalized Laplacian of network $G^{t}.$

From spectral graph theory \cite{Chung1996}, we know that due to the
connectedness of $A^{t},$ for every time point $t$ the solution to $L^{t}%
v_{1}^{t}=0$ corresponds to the first eigenvalue $\lambda_{1}^{t}=0$ which has
multiplicity one and the corresponding eigenvector $v_{1}$ is the eigenvector
$(D^{t})^{\frac{1}{2}}\boldsymbol{1}$, where $\boldsymbol{1}$ is the constant one vector and $D^{t}$ is the degree matrix for the adjacency matrix $A^{t}$.

Hence, the equation $\mathcal{L}^{0}v=0$ has a
$T-$dimensional subspace of solutions and we find its basis
explicitly:\ namely, for every $t$ we define the column vector $V^{t}%
\in\mathbb{R}^{NT},$ as a zero-padded vector with $v_{1}^{t}$ at the position
of the $t^{th}$ block. Thus, all solutions to $\mathcal{L}^{0}v=0$ are given by
$v=\sum_{t=1}^{T}\alpha_{t}V^{t}$ for arbitrary constants $\alpha_{t}$.

The main objective of the present paper is to consider an ideal case of a
temporal network which is slowly-changing in time, hence, is well approximated
by a temporal network following a \textit{constant block Jacobi model}: Let us consider the case where
$A^{t}=A$ for all $t$ and $W^{t,p}=W$ for all $t,p$. An important step in our construction is to "periodize" the temporal 
network, which will provide the existence of a nice closed-form solution of the resulting network. This is not a very 
artificial approach since the "slowly-changing" of the network assumes that the network does not vary too much from the initial to the final layer: Namely, we construct a "periodic"
supra-adjacency matrix $\mathcal{A}$ and its corresponding supra-Laplacian
matrix $\mathcal{L}$ for temporal networks, by including non-zero diagonal
blocks on the upper-right and lower-left corner blocks. In other words, we
include inter-layer weights between the first time layer $A^{1}$ and the last
time layer $A^{T}$. The resulting matrix $\mathcal{A}$ is a \textit{periodic constant block Jacobi} matrix which gives the name of the model. In view of the slowly-changing nature of the temporal
network $G^{t},$ the matrix $\mathcal{A}$ is a perturbation of the matrix
$\mathcal{A}^{0}$ and $\mathcal{L}$ is a perturbation of the matrix $\mathcal{L}%
^{0}.$

Furtheron, the resulting supra-Laplacian matrix $\mathcal{L}$ is given by the
following $T\times T$ block matrix, which may be easily proved to be an infinite periodic block
Jacobi matrix \cite{Sahbani2015}:
\begin{equation}
\mathcal{L}:=\underset{T}{\underbrace{\left(
\begin{array}
[c]{ccccc}%
\widetilde{L} & \widetilde{L}_{W} &  &  & \widetilde{L}_{W}\\
\widetilde{L}_{W} & \widetilde{L} & \widetilde{L}_{W} &  & \\
& \widetilde{L}_{W} & \widetilde{L} &  & \\
&  &  & \cdot\cdot\cdot & \widetilde{L}_{W}\\
\widetilde{L}_{W} &  &  & \widetilde{L}_{W} & \widetilde{L}%
\end{array}
\right)  }} \label{eq:ltildefinal}%
\end{equation}

We have to note that if we have the same
$\omega$ for all matrices $W,$ then the blocks of the block-diagonal matrix
$\mathcal{D}$ contain the matrices $D^{t}+2\omega I$. Since for every $t$
holds equation $L^{t}=I-D^{-1/2}AD^{-1/2},$ and since the matrix $D^{-1/2}AD^{-1/2}$ has entries
$d_{i}^{-1/2}d_{j}^{-1/2}a_{ij}$, we see that $\widetilde{L}$ is a
perturbation of $L$ which has just the elements $-\left(  d_{i}+2\omega
\right)  ^{-1/2}\left(  d_{j}+2\omega\right)  ^{-1/2}a_{ij}$ and not
$-d_{i}^{-1/2}d_{j}^{-1/2}a_{ij}$. Hence, written formally, we have the equality
\[
\widetilde{L}=I-\left(  D+2\omega I\right)  ^{-1/2}A\left(  D+2\omega
I\right)  ^{-1/2}%
\]
On the other hand, the matrix $\widetilde{L}_{W}$ is equal to $-\omega\left(  D+2\omega
I\right)  ^{-1}$ , in equation (\ref{eq:ltildefinal}).

The big advantage of the constant block Jacobi model is that we can find "explicitly" its 
spectrum which we discuss in the next sections. 

\section{Smallest eigenvalues and paired eigenvectors of the supra-Laplacian $\mathcal{L}$ of temporal networks following constant block Jacobi model\label{sec3}}
As we know from spectral graph theory \cite{Chung1996}, the eigenvalues of
the Laplacian $L^{t}$ and of the supra-Laplacian $\mathcal{L}$ are
non-negative, and the minimal eigenvalue is $0$, as mentioned above. As usual,
in the applications the small eigenvalues and the corresponding eigenvectors
are of particular importance. By perturbation theory, some of those
eigenvalues which are very close to $0$ are obtained as a direct perturbation
of the $0$ eigenvalues of all separate time layer Laplacian matrices $L^{t},$
and the same holds about their paired eigenvectors. On the other hand,
the eigenvectors paired to the bigger eigenvalues are obtained as
perturbations not only of the $0$ eigenvectors of the separate matrices
$L^{t}$ but also of the Fielder (and the higher) eigenvectors of the separate
matrices $L^{t}$.

The solution for the Laplacian $\mathcal{L}$ in equation (\ref{eq:ltildefinal})
is defined by:
\begin{equation}
\mathcal{L}\psi=\lambda\psi\label{ABreduced}%
\end{equation}
and for finding it we apply a classical technique based on discrete Fourier transforms
(DFTs), see e.g. \cite{Sahbani2015}. To do this we
represent each vector $\psi\in\mathbb{R}^{NT}$ as the sequence of vectors
$\left[  \psi_{1},\psi_{2},...,\psi_{T}\right]  $ where each vector $\psi_{j}$
is the portion of eigenvector $\psi$ corresponding to the $j^{th}$ time block.
Then equation (\ref{ABreduced}) splits into the equations
\begin{equation}
\widetilde{L}_{W}\psi_{j-1}+\widetilde{L}\psi_{j}+\widetilde{L}_{W}\psi
_{j+1}=\lambda\psi_{j}\qquad\text{for }j=1,2,...,T \label{ABreducedDETAIL}%
\end{equation}
where for the sake of notation simplicity we have put
\[
\psi_{0}=\psi_{T},\qquad\psi_{T+1}=\psi_{1}.
\]
For $k=0,1,2,...,T-1,$ we denote the DFT of
vector $\psi$ at value $k$ by $\widehat{\psi}(k)\in\mathbb{R}^{N},$ and put
\begin{equation}
\widehat{\psi}(k):=\sum_{j=0}^{T-1}%
e^{-ijk\frac{2\pi}{T}}\psi_{j+1}. \label{ABreducedE}%
\end{equation}
It is important that from the set of DFT vectors $\{\widehat{\psi}\left(k\right)\}_{k=0}^{T-1}$ we may recover the whole vector
$\psi\in\mathbb{R}^{NT}$ using the Fourier inversion formula:
\begin{equation}
\psi_{j}=\frac{1}{T}\sum_{k=0}^{T-1}\widehat{\psi}(k)e^{ijk\frac{2\pi}{T}}. \label{DFT-Inversion}%
\end{equation}

Now by applying the DFT (\ref{ABreducedE}) to equations (\ref{ABreducedDETAIL}%
) (i.e. by multiplying by exponents and summing up the equations), we obtain
the fundamental equations satisfied by the DFT of the vector $\psi$ defined in
formula (\ref{ABreducedE}):
\begin{eqnarray}
\left[  \widetilde{L}+2\cos\left(  k\frac{2\pi}{T}\right)  \widetilde{L}%
_{W}\right]  \widehat{\psi}\left(  k\right)  =\lambda\widehat{\psi}\left(
k\right)\hspace{0.25em} & \label{B+coskA} \\ \text{for }& k=0,1,...,T-1. \nonumber 
\end{eqnarray}

The following theorem justifies the application of the DFTs for solving the
system (\ref{ABreduced}):

\begin{theorem}
\label{Tsahbani} The spectrum (with multiplicities) of the supra-Laplacian
$\mathcal{L}$ in equation (\ref{eq:ltildefinal}) of a temporal network following a periodic constant block Jacobi model
coincides with the union of the spectra of the matrices $\widetilde{L}%
+2\cos\left(  k\frac{2\pi}{T}\right)  \widetilde{L}_{W},$ i.e.
\begin{equation}
spec\left(  \mathcal{L}\right)  =\cup_{k=0}^{T-1}spec\left(  \widetilde
{L}+2\cos\left(  k\frac{2\pi}{T}\right)  \widetilde{L}_{W}\right)
\label{SpecUnion}%
\end{equation}

\begin{proof}
First, we prove the inclusion
\[
spec\left(  \mathcal{L}\right)  \subseteq\cup_{k=0}^{T-1}spec\left(
\widetilde{L}+2\cos\left(  k\frac{2\pi}{T}\right)  \widetilde{L}_{W}\right)  .
\]
Indeed, by the above arguments, if we have an eigenvalue $\lambda$ with
eigenvector $\psi$ solving system (\ref{ABreducedDETAIL}), then for every $k$
with $0\leq k\leq T-1$ we have equation (\ref{B+coskA}), i.e.
\[
\left[  \widetilde{L}+2\cos\left(  k\frac{2\pi}{T}\right)  \widetilde{L}%
_{W}\right]  \widehat{\psi}\left(  k\right)  =\lambda\widehat{\psi}\left(
k\right)  .
\]
Hence, $\lambda$ is an eigenvalue for all matrices $\widetilde{L}+2\cos\left(
k\frac{2\pi}{T}\right)  \widetilde{L}_{W}$ with eigenvector $\widehat{\psi
}\left(  k\right).$ Now, we prove the opposite inclusion:
\[
\cup_{k=0}^{T-1}spec\left(  \widetilde{L}+2\cos\left(  k\frac{2\pi}{T}\right)
\widetilde{L}_{W}\right)  \subseteq spec\left(  \mathcal{L} \right)  .
\]
Assume that $\lambda^{\ast}$ is an eigenvalue with eigenvector $v^{\ast}$ for
the matrix $\widetilde{L}+2\cos\left(  k\frac{2\pi}{T}\right)  \widetilde
{L}_{W}$, i.e.
\[
\left[  \widetilde{L}+2\cos\left(  k\frac{2\pi}{T}\right)  \widetilde{L}%
_{W}\right]  v^{\ast}=\lambda^{\ast}v^{\ast}.
\]
We define the vector $\varphi\in\mathbb{R}^{NT}$ by putting
\begin{align*}
\varphi_{k+1} &  =v^{\ast}\\
\varphi_{m} &  =0\qquad\text{for }m\neq k+1,m=1,2,...,T.
\end{align*}
By the inversion formula (\ref{DFT-Inversion}) we define the vector
\[
\psi_{j}:=\varphi_{k+1}e^{ijk\frac{2\pi}{T}}\qquad\text{for }j=1,2,...,T.
\]
We show that it satisfies the eigenvalue equation (\ref{ABreducedDETAIL})
since
\[
\widetilde{L}_{W}\psi_{j-1}+\widetilde{L}\psi_{j}+\widetilde{L}_{W}\psi
_{j+1}=\lambda^{\ast}\psi_{j}%
\]
i.e.
\[
e^{i\left(  j-1\right)  k\frac{2\pi}{T}}\widetilde{L}_{W}v^{\ast}%
+e^{ijk\frac{2\pi}{T}}\widetilde{L}v^{\ast}+e^{i\left(  j+1\right)
k\frac{2\pi}{T}}\widetilde{L}_{W}v^{\ast}=\lambda^{\ast}e^{ijk\frac{2\pi}{T}%
}v^{\ast}%
\]
But the last is equivalent to equation
\[
e^{-ik\frac{2\pi}{T}}\widetilde{L}_{W}v^{\ast}+\widetilde{L}v^{\ast
}+e^{ik\frac{2\pi}{T}}\widetilde{L}_{W}v^{\ast}=\lambda^{\ast}v^{\ast}%
\]
hence, to equation $\widetilde{L}v^{\ast}+2\cos\left(  k\frac{2\pi}{T}\right)
\widetilde{L}_{W}v^{\ast}=\lambda^{\ast}v^{\ast};$ which was our assumption.
This completes the proof. 
\end{proof}
\end{theorem}

\begin{figure}
\includegraphics[width=.99\columnwidth]{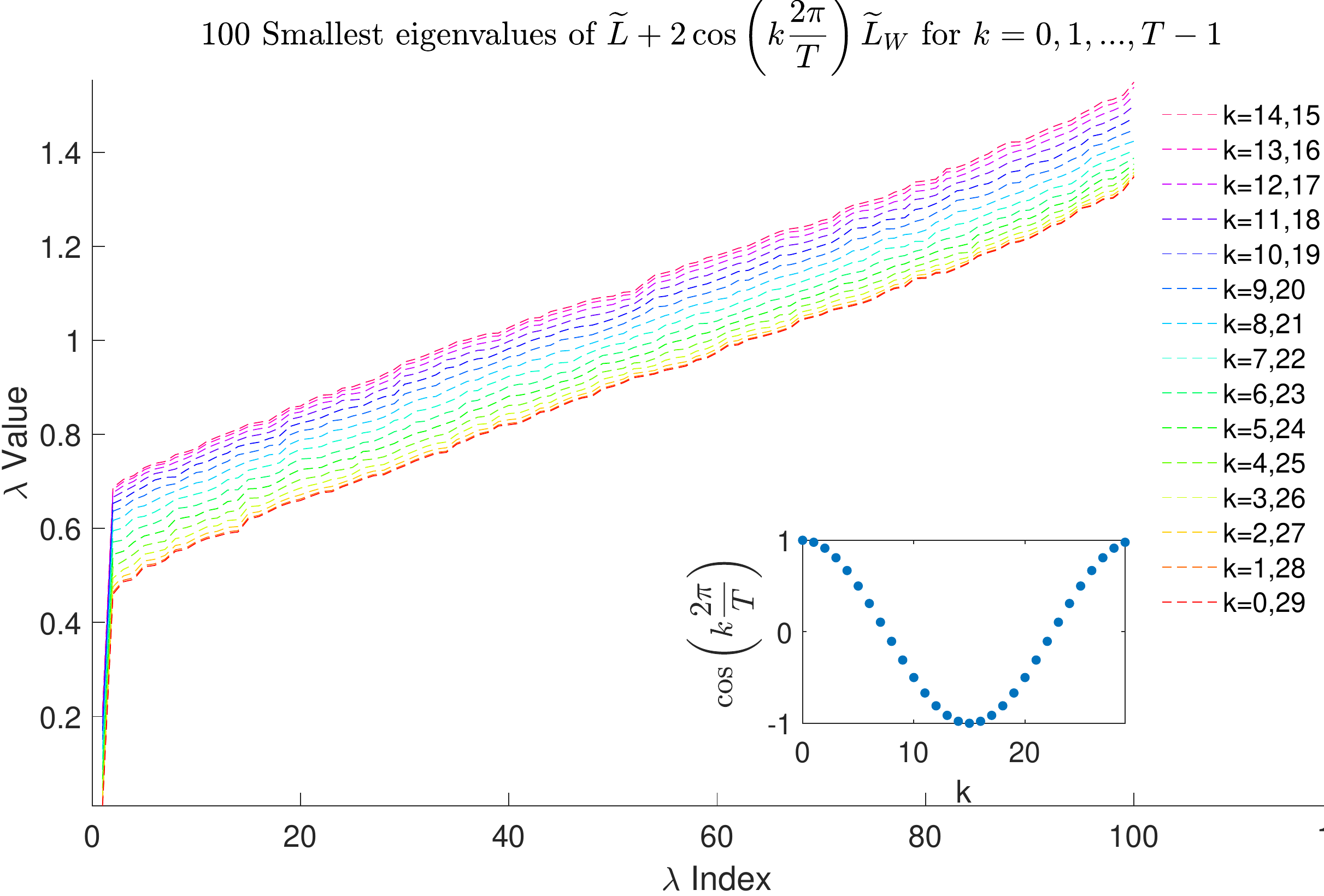}
\caption{\label{fig:eigenest_ER}\textbf{The \pmb{$100$} smallest eigenvalues of matrices \pmb{$\tilde{L}+2\cos\left(
k\frac{2\pi}{T}\right)  \tilde{L}_{W}$} for each \pmb{$k=0,1,2,...,29$.}} The matrices
$\tilde{L}$ and $\tilde{L}_{W}$ are obtained from a temporal benchmark network composed
of $T=30$ \textbf{Erdos-Renyi random graphs} each with $N=100$ nodes and edge probability $p=0.3$ (such dense consecutive ER networks are slowly-changing). The inter-layer weights $\omega$ are fixed at $1$. We include the
additional plot of $\cos\left(  k\frac{2\pi}{T}\right)  $ which determines the
monotonically increasing behavior of eigenvalues corresponding to $0\leq k\leq14$ and
monotonically decreasing behaviour of eigenvalues corresponding to $15\leq k\leq29$.}
\end{figure}

In Figure \ref{fig:eigenest_ER} we have displayed the first $100$ eigenvalues of
the matrix $L=\widetilde{L}+2\cos\left(  k\frac{2\pi}{T}\right)  \widetilde
{L}_{W}$ from equation (\ref{B+coskA}), where we see that for every $j\geq1,$
the $j^{th}$ eigenvalue $\lambda_{j}^{\left(  k\right)  }$ of all matrices
$\widetilde{L}+2\cos\left(  k\frac{2\pi}{T}\right)  \widetilde{L}_{W}$ is
monotonically increasing with $k$ for
\begin{align*}
0 &  \leq k\leq\frac{T-1}{2}-1\text{ if }T\text{ is odd}\\
&  \text{and}\\
0 &  \leq k\leq\frac{T}{2}-1\text{ if }T\text{ is even}.
\end{align*}
The following proposition explains the behavior of the eigenvalues.

\begin{proposition}
Without loss of generality assume that $T$ is odd. Then the $j^{th}$
eigenvalues of the matrices $\widetilde{L}+2\cos\left(  k\frac{2\pi}%
{T}\right)  \widetilde{L}_{W}$ satisfy
\[
\lambda_{j}^{\left(  0\right)  }\leq\lambda_{j}^{\left(  1\right)  }\leq
\cdot\cdot\cdot\leq\lambda_{j}^{\left(  \frac{T-1}{2}-1\right)  }.
\]

\begin{proof}
The proof of this proposition is direct consequence of Theorem 8.1.5. in
\cite{Golub1996} which states that for symmetric matrices $V$ and $E$ of size
$N\times N,$ and for all eigenvalues $\lambda_{j}$, for $j=1,2,...,N,$
hold the inequalities:
\begin{equation}
\lambda_{j}\left(  V\right)  +\lambda_{\min}\left(  E\right)  \leq\lambda
_{j}\left(  V+E\right)  \leq\lambda_{j}\left(  V\right)  +\lambda_{\max
}\left(  E\right)  .\label{lambdakP}%
\end{equation}

We take into account the fact that the eigenvalues of the diagonal matrix $\widetilde{L}%
_{W}$ are non-negative since they coincide with all non-negative weights
$\omega_{j}^{t,p}$. In particular, if all they are equal to a constant
$\omega$, then we see that
\[
\lambda_{j}^{k}=\lambda_{j}(\widetilde{L})+2\cos\left(  k\frac{2\pi}%
{T}\right)  \omega.
\]
This completes the proof. 
\end{proof}
\end{proposition}

\begin{figure}
\centering
\includegraphics[width=.49\columnwidth]{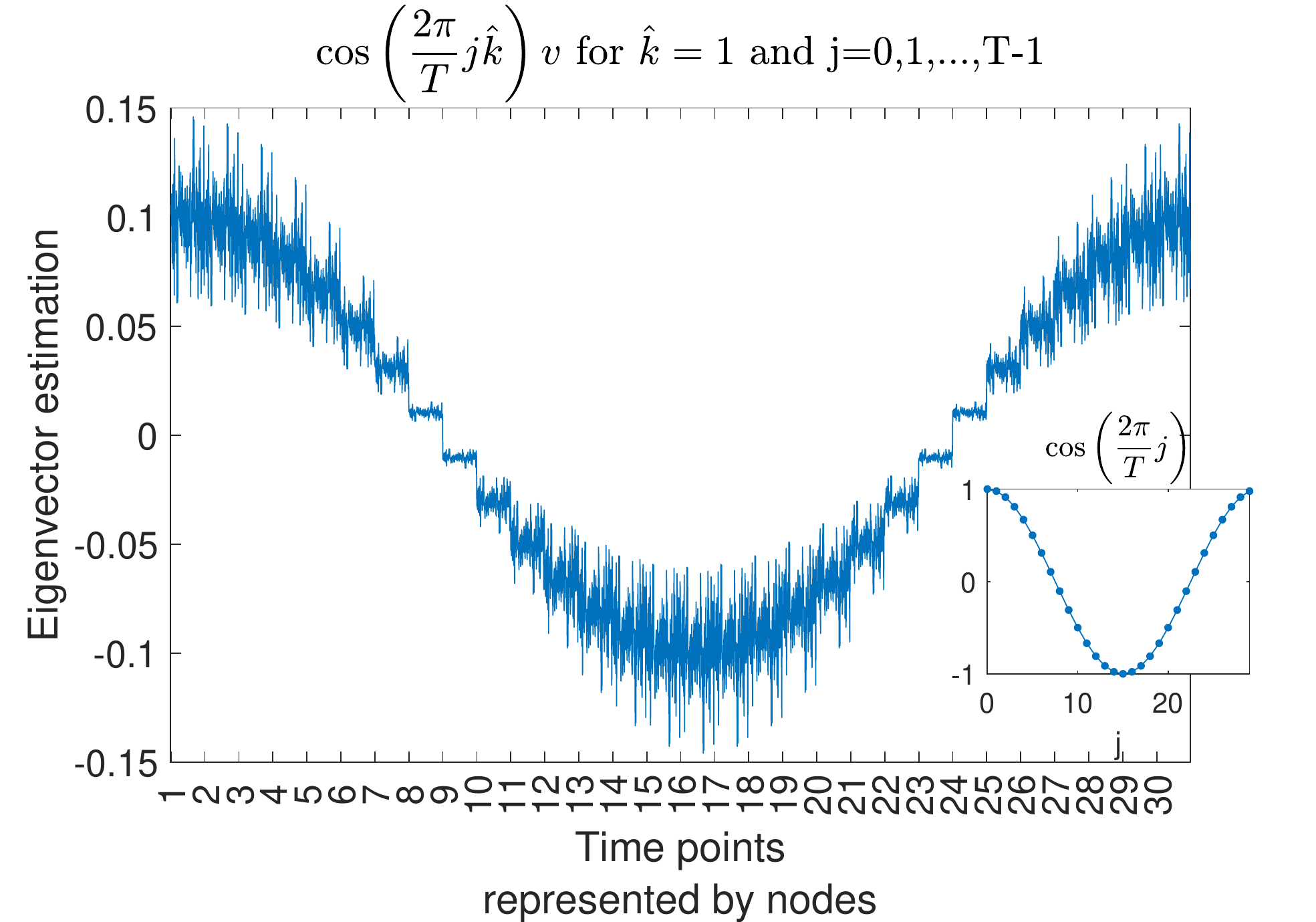}\hfill
\includegraphics[width=.49\columnwidth]{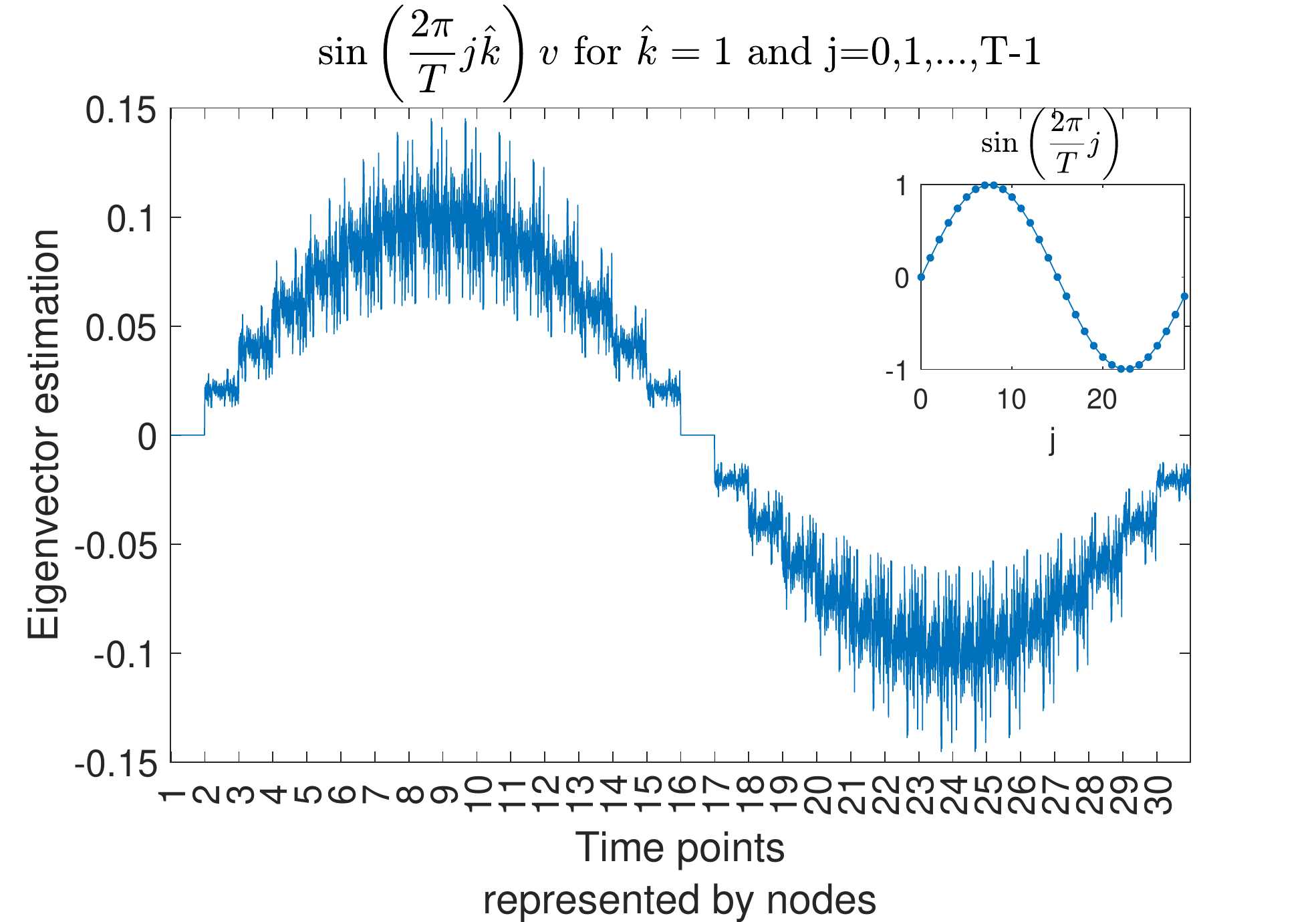}\hfill 
\includegraphics[width=.49\columnwidth]{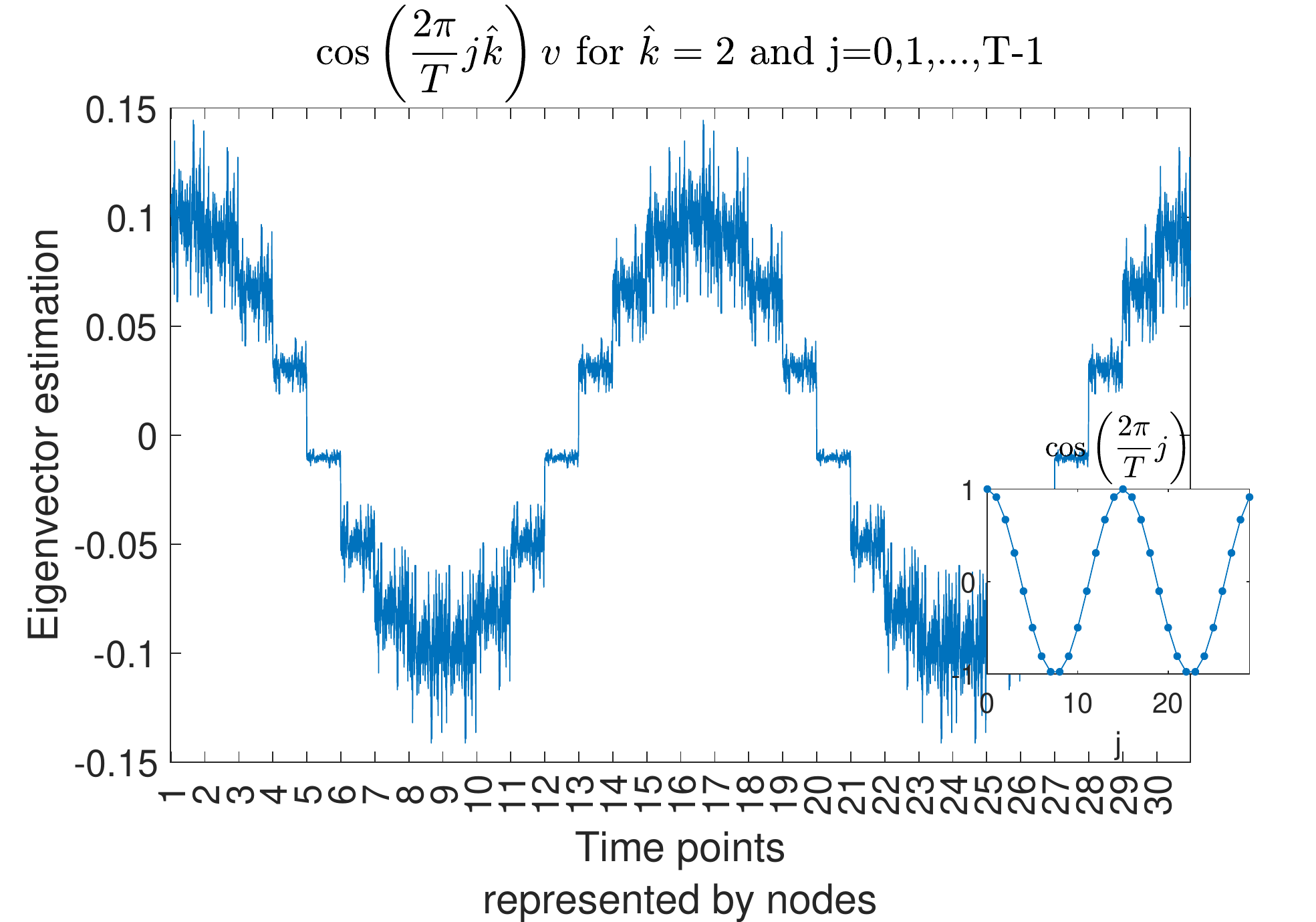}\hfill 
\includegraphics[width=.49\columnwidth]{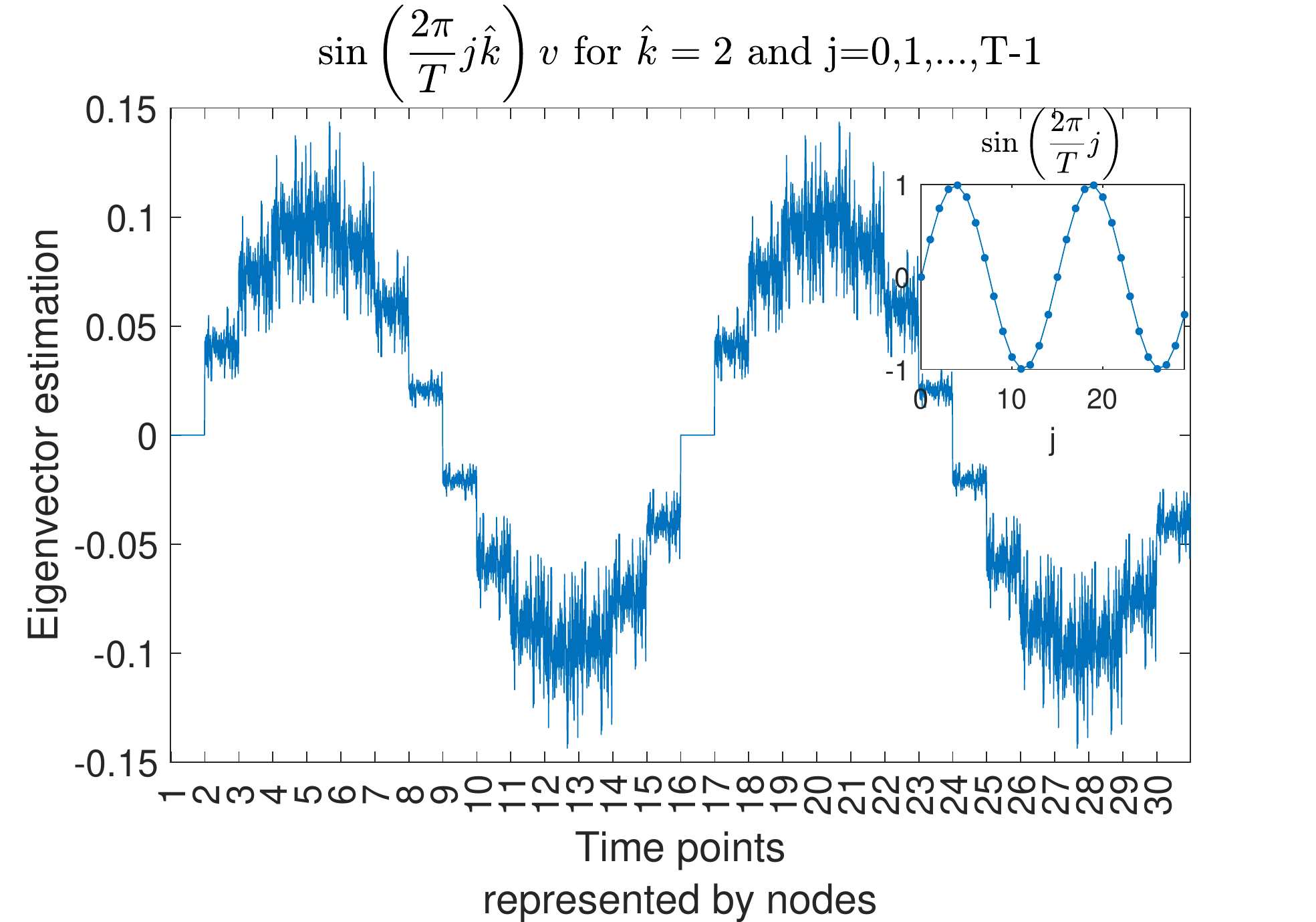}\hfill 
\includegraphics[width=.49\columnwidth]{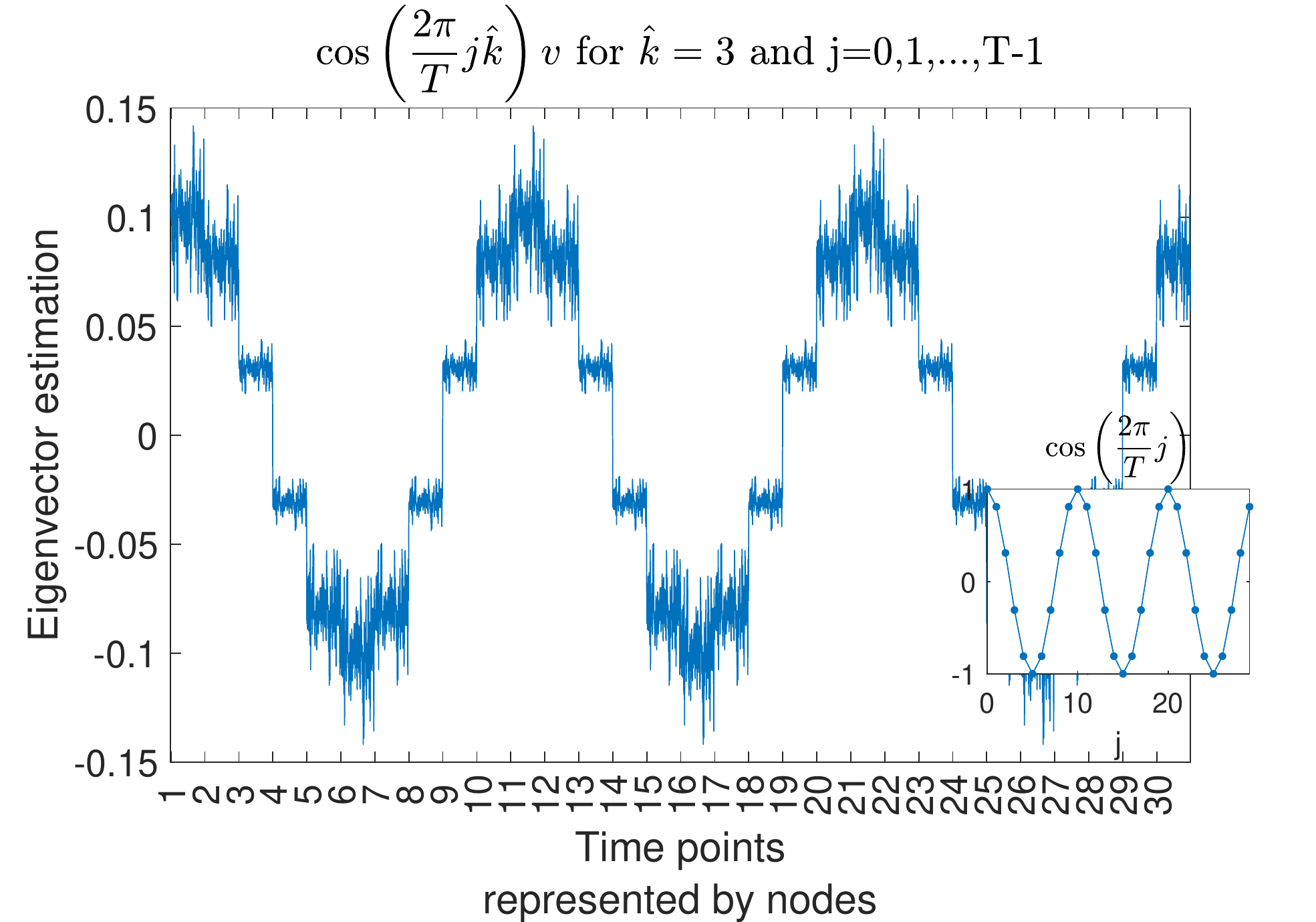}\hfill 
\includegraphics[width=.49\columnwidth]{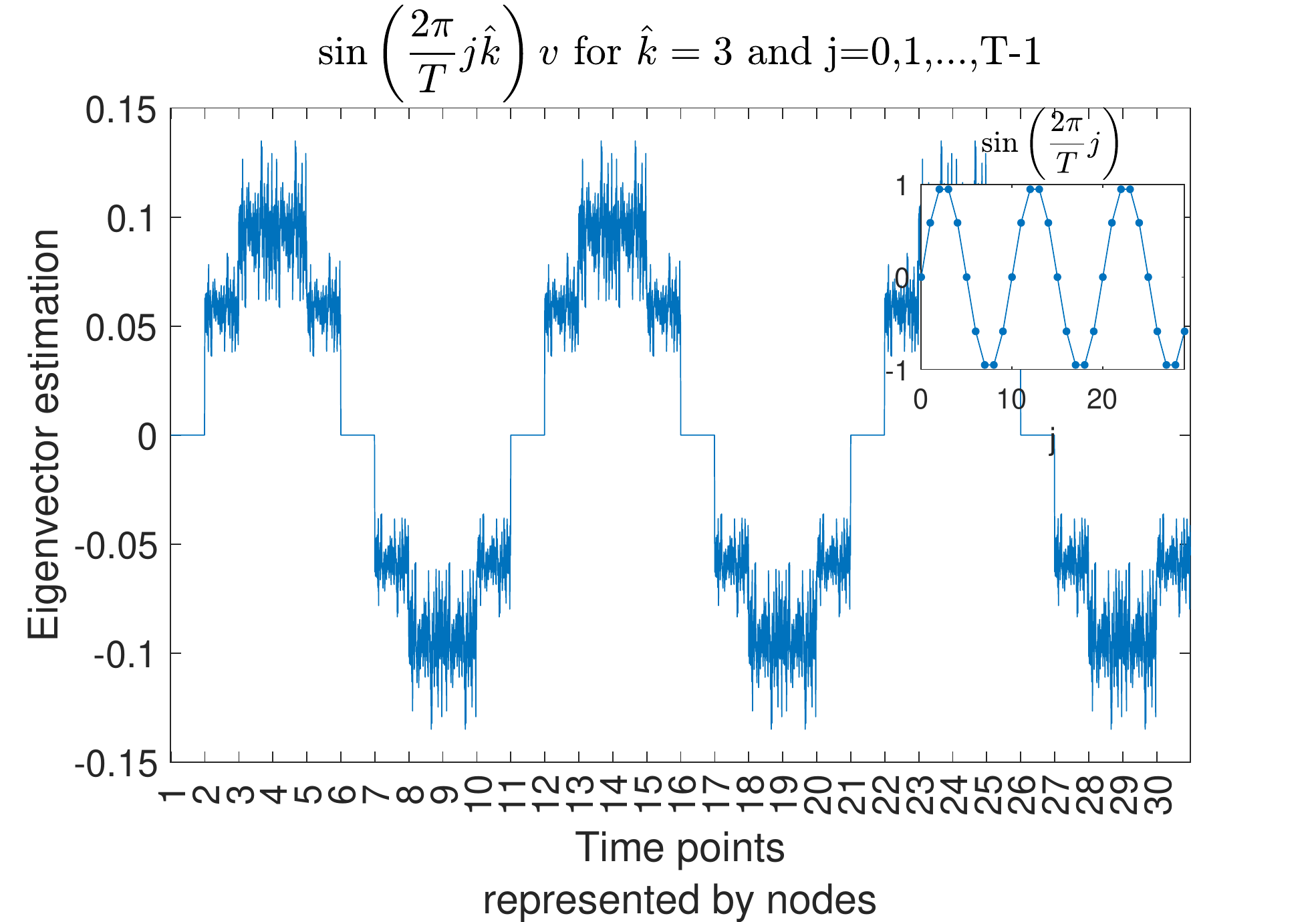}
\caption{\label{fig:CosVsSin}\textbf{Eigenvector estimations for supra-Laplacian matrix \pmb{$\mathcal{L}$}.} This figure visualizes eigenvectors from
equation (\ref{psijk1}) for $\hat{k}=1,2,3,$ each accompanied by the corresponding
graph of the $\cos$ and $\sin$ functions. The eigenvector $v$ corresponds to
the eigenvalue $\lambda=0$ which is a solution to the eigenvalue problem
(\ref{B+coskA}). The matrices $\tilde{L}$ and $\tilde{L}_{W}$ are obtained
from a temporal network following the constant block Jacobi model composed
of $T=30$ \textbf{Erdos-Renyi random graphs} each with $N=100$ nodes and edge probability $p=0.3$. The inter-layer weights $\omega$ are fixed at $1$.}%
\end{figure}

\begin{figure}[t]
\centering
\includegraphics[width=.99\columnwidth]{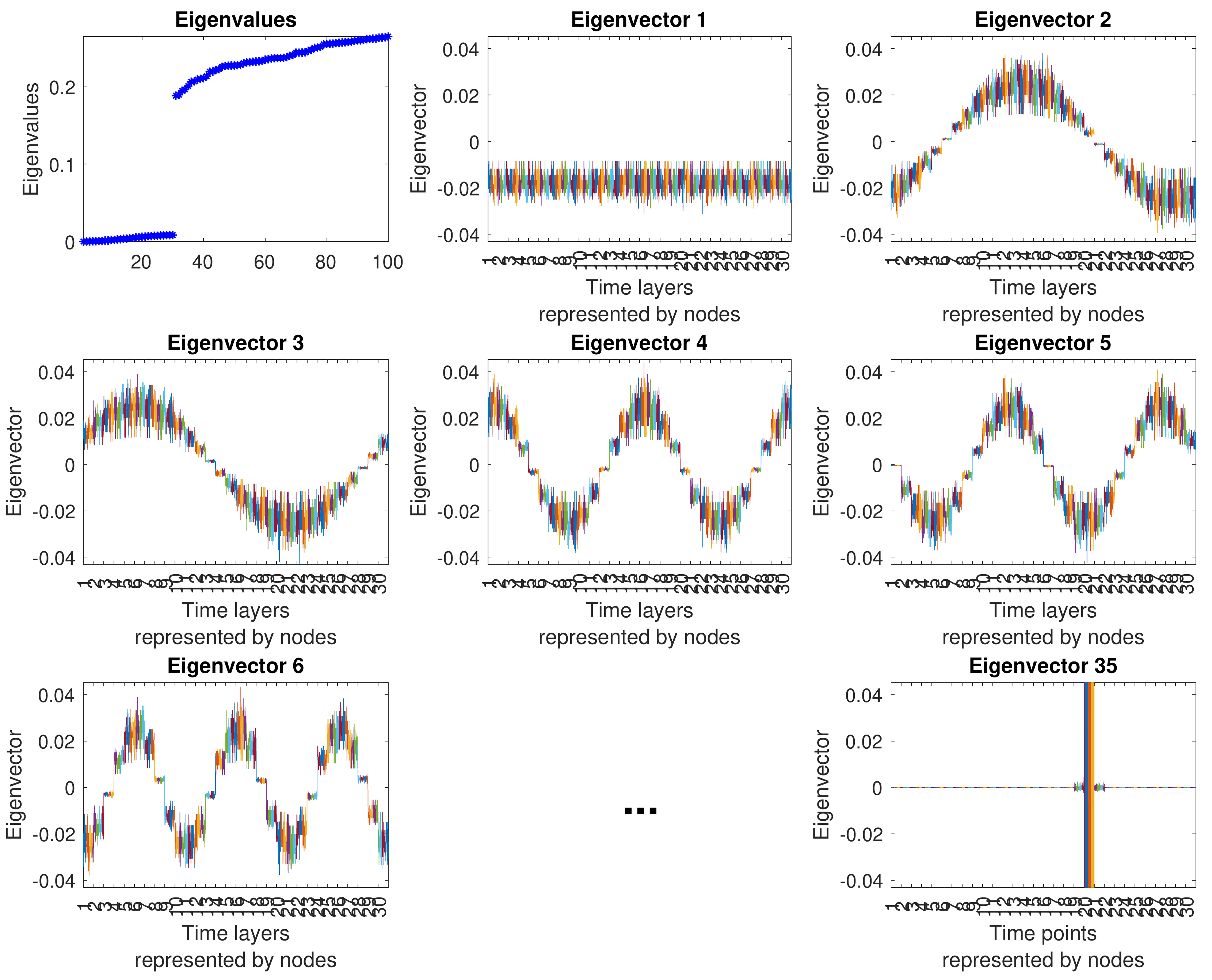}\caption{\label{fig:ER_eigen} \textbf{Eigenvalues and
eigenvectors for an \textbf{Erdos-Renyi} benchmark temporal network.} The Erdos-Renyi temporal benchmark
network is composed of $T=30$ random Erdos-Renyi graphs with $N=100$ nodes and
$p=0.1$ edge probability. The inter-layer weights are set to $\omega=0.01$. We
plot the $100$ smallest eigenvalues of the corresponding supra-Laplacian
matrix, the $6$ eigenvectors corresponding to the $6$ smallest eigenvalues and
the $35^{th}$ eigenvector. The jump of the eigenvalue graph indicates precisely the position of $\lambda^{*}$ for index 31 and all following eigenvectors look as the $35^{th}$ eigenvector plotted which captures local variability.  Colouring of each eigenvector is consistent with
the components that belong to different time points.}
\end{figure}

Now, by means of Theorem~\ref{Tsahbani}, we show how to construct a solution
to eigenvalue equation (\ref{ABreduced}) by using equality (\ref{B+coskA}): Fix a $k=\hat{k}$ and consider an
eigenvector $v$ with eigenvalue $\hat{\lambda}$ solving the eigenvalue problem
(\ref{B+coskA}) for $k=\hat{k}.$ We assume that $\hat{\lambda}$ is among the
smallest eigenvalues, close to $0.$ We are seeking for a block-vector
$\Psi=\left(  \psi_{1},\psi_{2},...,\psi_{T}\right)  \in\mathbb{R}^{NT}$ for
which $\widehat{\Psi}\left(  k\right)  =\varphi_{k},$ where the block-vector
$\Phi=\left(  \varphi_{1},...,\varphi_{T}\right)  \in\mathbb{R}^{NT}$ is
defined as
\[
\varphi_{k}:=\left\{
\begin{array}
[c]{cc}%
v & \quad\text{for }k=\hat{k}\\
0 & \quad\text{for }k\neq\hat{k}%
\end{array}
\right.
\]
Now we apply the inversion formula (\ref{DFT-Inversion}) to the vector $\Phi,$
and obtain the block-vector $\Psi\in\mathbb{C}^{NT}$ with components
\begin{equation}
\psi_{j}=e^{\frac{2\pi}{T}ij\hat{k}}v\qquad\text{for }%
j=0,1,...,T-1.\label{psij}%
\end{equation}
Thus we have $\varphi_{k}=0$ for $k\neq\hat{k},$ and $\Psi$ is a solution to
the eigenvalue equation (\ref{ABreduced}) with the same $\hat{\lambda}.$ Since
the vector $\Psi$ is complex valued, we obtain two real-valued vectors
($\in\mathbb{R}^{NT}$), by taking the real and imaginary parts of
$e^{\frac{2\pi}{T}ij\hat{k}},$ namely:
\begin{align}
\psi_{j}^{R} &  :=\cos\left(  \frac{2\pi}{T}j\hat{k}\right)  \times
v\qquad\text{for }j=0,1,...,T-1\label{psijk1}\\
\psi_{j}^{I} &  :=\sin\left(  \frac{2\pi}{T}j\hat{k}\right)  \times
v\qquad\text{for }j=0,1,...,T-1\nonumber
\end{align}

In Figure~\ref{fig:CosVsSin} we visualise solutions (\ref{psijk1}) for
$\hat{k}=1,2,3,$ accompanied by the corresponding plots of $\cos(\frac{2\pi}
{T}j\hat{k})$ and $\sin(\frac{2\pi}{T}j\hat{k})$ for $j=0,1,...,T-1$.

Every eigenvalue in equation (\ref{B+coskA}) has \textit{even} multiplicity
due to the equality of the two matrices as indicated below:
\begin{align}
\widetilde{L}+2\cos\left(  k\frac{2\pi}{T}\right)  \widetilde{L}_{W}= &\widetilde{L}+2\cos\left(  \left(  T-k\right)  \frac{2\pi}{T}\right)
\widetilde{L}_{W}\nonumber\\&\text{for }0\leq k\leq\frac{T-1}{2}-1\nonumber;
\end{align}
the double multiplicity of the eigenvalues is clearly observed in Figure
\ref{fig:eigenest_ER}. In the case of odd $T$ there are unique eigenvalues just
for $k=\frac{T-1}{2}-1$; for even $T$ all eigenvalues have even multiplicity.
For $\hat{k}=0$ we have one solution $\Psi$ with $\psi_{j}=v$ corresponding to
the zero eigenvalue, $\hat{\lambda}=0.$

By using the results of perturbation theory for invariant subspaces
\cite{Golub1996,Luxburg} we see that for every eigenvalue with even
multiplicity, we may estimate the perturbation of its eigenspace, i.e. the
space of its eigenvectors. Thus we obtain the solutions which look like
\textquotedblleft block sinusoids\textquotedblright\ of $\cos$ and $\sin$
type, Figure~\ref{fig:CosVsSin}. The perturbation of the two-dimensional space
spanned by $\cos$ and $\sin$ type solutions, results in a two-dimensional
space corresponding to the perturbed eigenvalue of the matrix $\mathcal{L}$.
These eigenvectors may differ from $\cos$ or $\sin$ type solutions. 

The above theoretical results have a direct impact on the eigenvectors of the supra-Laplacian $\mathcal{L}$, Figure~\ref{fig:ER_eigen}.
We show that the eigenvectors corresponding to the
eigenvalues of the supra-Laplacian $\mathcal{L}$, which are close to $0,$ are
obtained by perturbation of the eigenvectors corresponding to the $0$
eigenvalues of the separate layers $L^{t}$, derived as $\left(
D^{t}\right)  ^{\frac{1}{2}}\boldsymbol{1}$. Thus they do not carry any
information about the finer description of that layer as does the
Fiedler vector. These eigenvectors of $\mathcal{L}$ give us only information about
all $T$ time layers being separate from each other. The
bigger eigenvalues of $\mathcal{L}$ have eigenvectors which are perturbations of
mixtures of higher eigenvectors for networks $L^{t}$, i.e. they contain
information from the Fiedler eigenvectors for the separate networks $L^{t}$.
We can conclude that only after the block nature of the constant block Jacobi model 
in the temporal network is captured the eigenvectors start capturing variability introduced 
by some certain within-layer patterns, which is clearly seen from Figure~\ref{fig:ER_eigen}.

\begin{figure}
\includegraphics[width=.49\columnwidth]{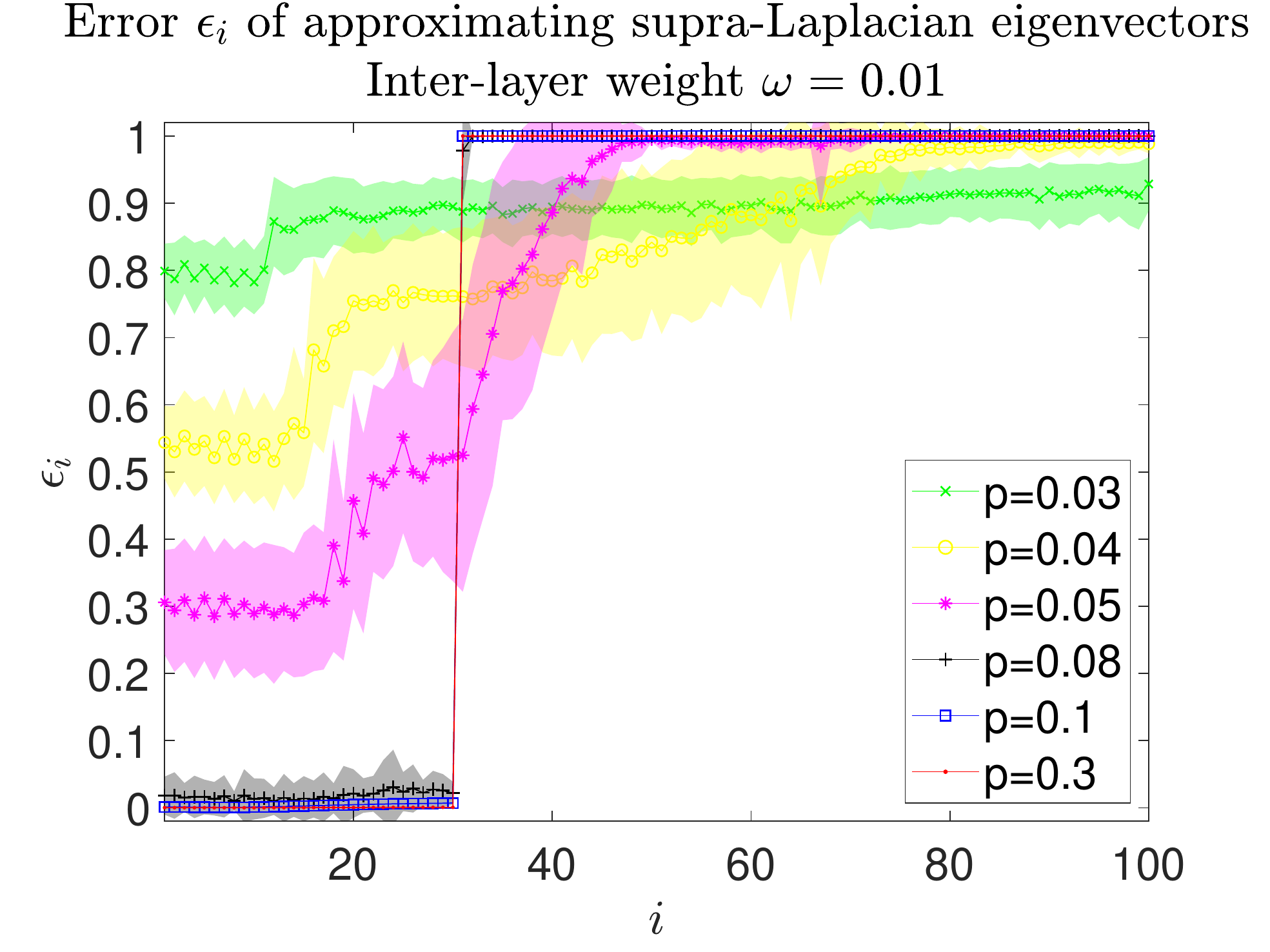}\hfill
\includegraphics[width=.49\columnwidth]{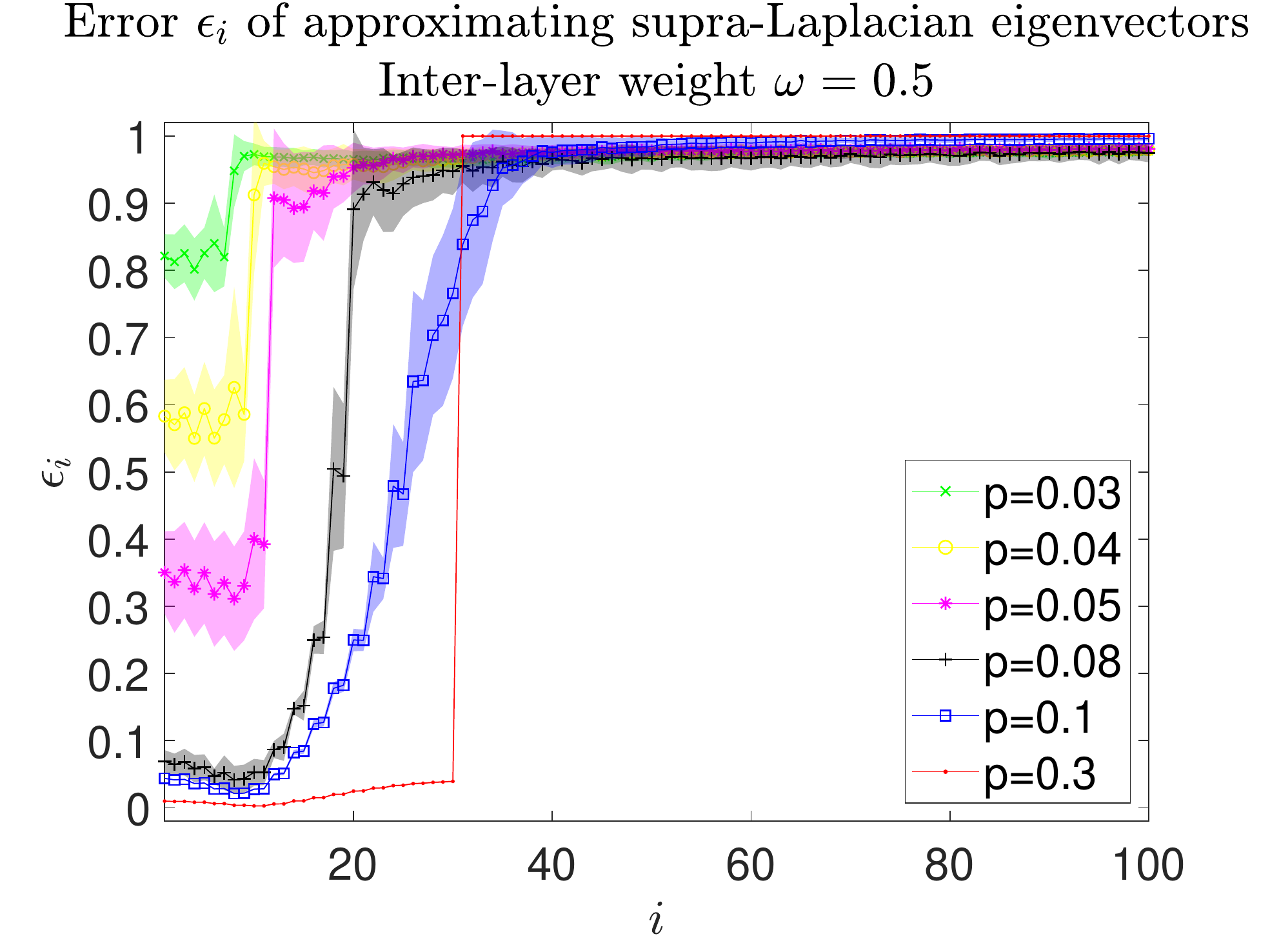}\hfill
\includegraphics[width=.49\columnwidth]{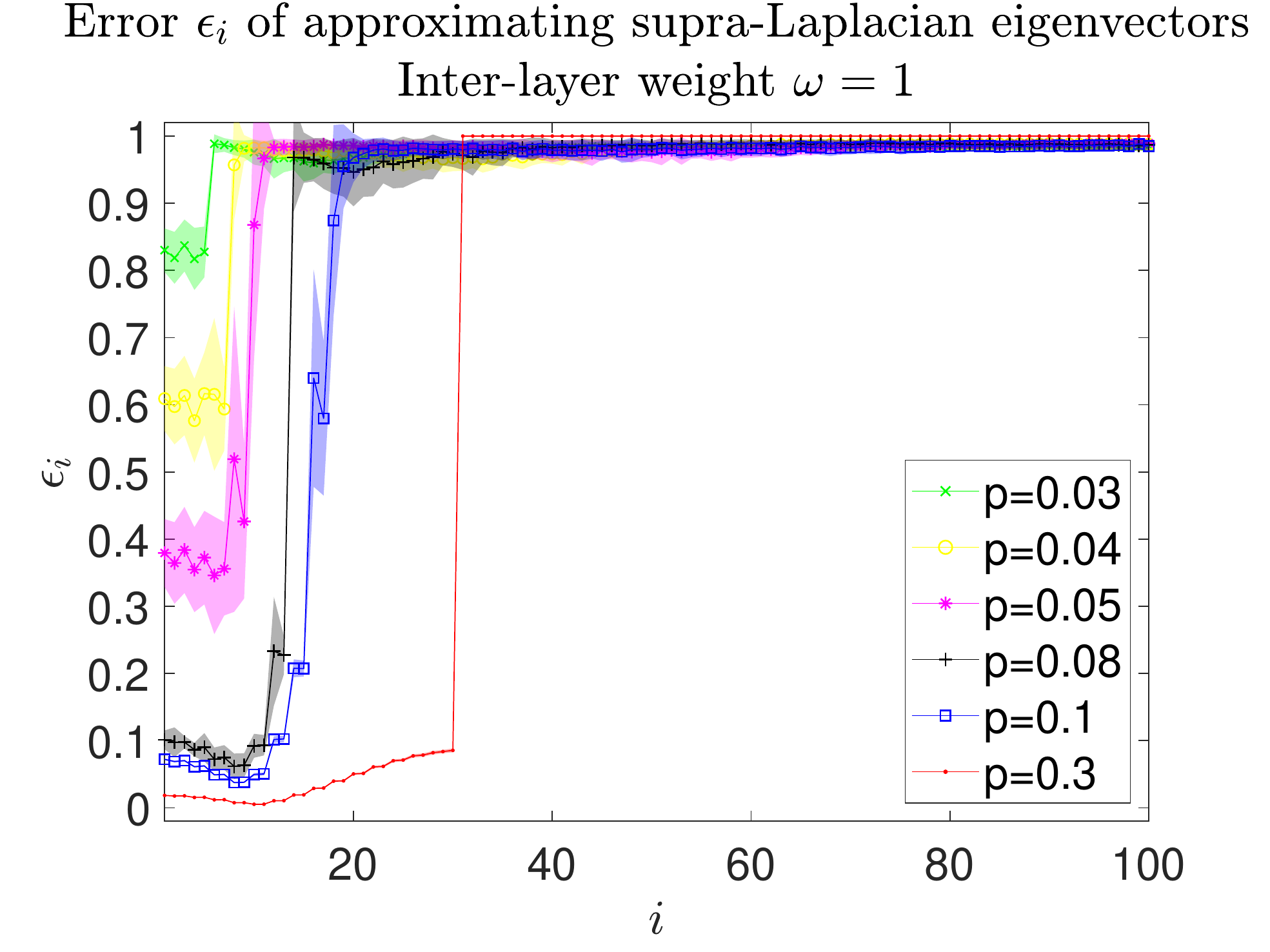}\hfill
\includegraphics[width=.49\columnwidth]{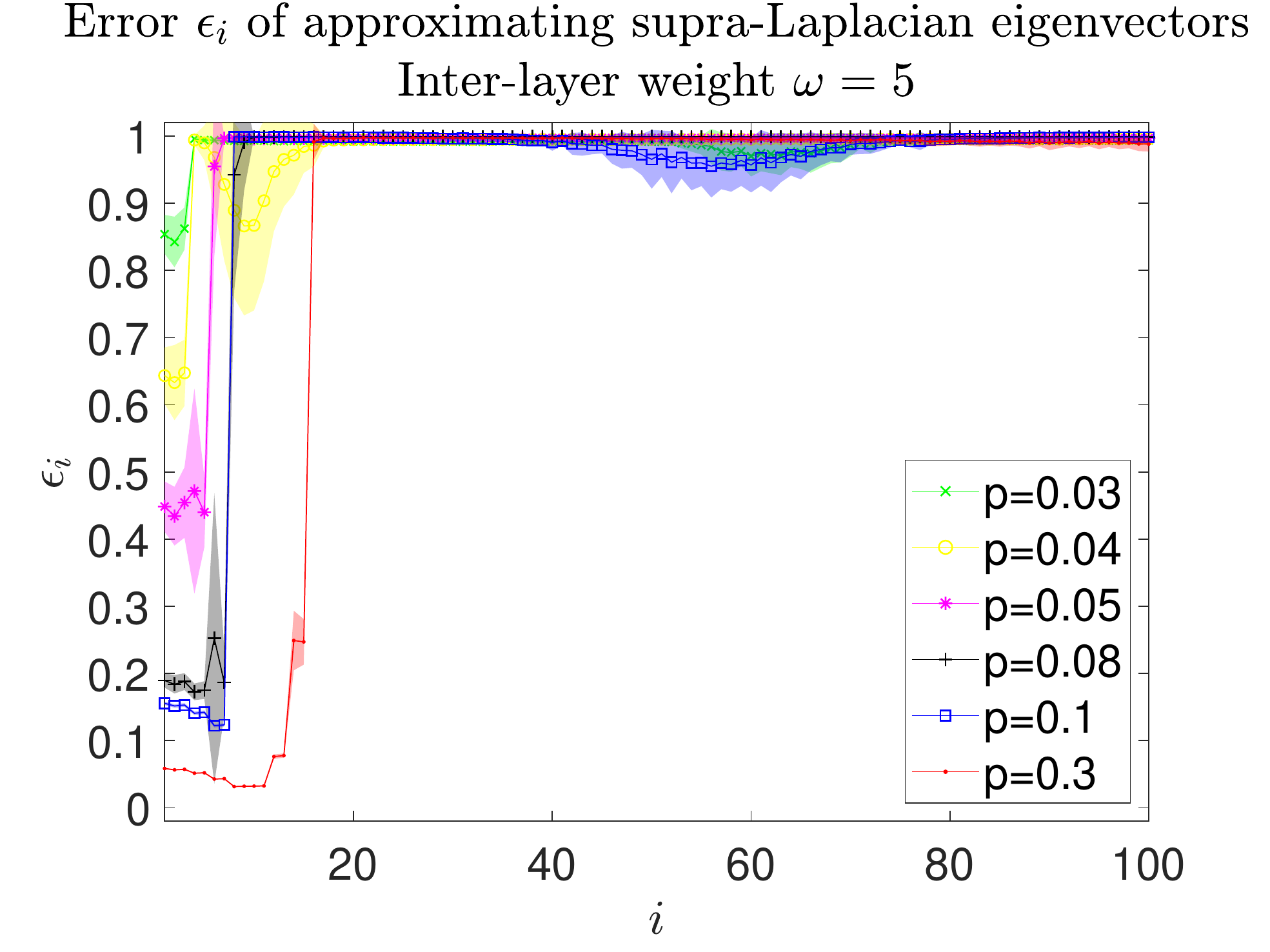}\hfill
\caption{\label{fig:err}\textbf{Error $\epsilon_{i}$
of approximating supra-Laplacian eigenvectors (corresponding to eigenvalue
\pmb{$\lambda_{i}$} for \pmb{$i=1,2,3,....,TN$)} by their separate time layers eigenvectors for
the benchmark temporal network.} All of the benchmark temporal networks were simulated using $T=30$ \textbf{random Erdos-Renyi graphs} with $N=100$ nodes and varying edge probabilities $p=0.03,0.04,0.05,0.08,0.1,0.3$ edge probabilities. Each of the four plots captures the results for different inter-layer weights set to $\omega=0.01,0.05,1,5$. For each parameter combination $(p,\omega)$ we simulate $100$ networks and show their average error $\epsilon_{i}$ with $1$ st.dev. intervals. The obtained approximation average errors and st.dev. intervals are visualized for the first $100$ eigenvectors although at most $T+1$ regressions are needed to capture all $T$ layers as separate layers.} 
\end{figure}
\section{Properties of the eigenvectors corresponding to small eigenvalues of the supra-Laplacian $\mathcal{L}$\label{sec5}}

In this section we empirically showcase the theoretical results that eigenvectors corresponding 
to the small eigenvalues of $\mathcal{L}$ are well-approximated by linear combinations of the 
eigenvectors (paired to the zero eigenvalue) of the separate layers. We investigate their behavior 
with respect to the edge density of the layers and the inter-layer weights. 

\subsection{Evaluating the approximation of the eigenvectors of $\mathcal{L}$ using the eigenvectors of the separate time layers}
Let $\overline{\Lambda}$ be the set of smallest eigenvalues with paired
eigenvectors well-approximated by the subspace of eigenvectors
corresponding to the $0$ eigenvalues for the separate layers. The
theoretical results from Sec.~\ref{sec3} guarantee that the
eigenvectors $v$ corresponding to $\lambda\in\overline{\Lambda}$ satisfy (see
Sec.~\ref{sec2} for $V^{t}$ def.)%
\begin{equation}
\min_{\left\{  \alpha_{t}\right\}  }\left\Vert v-\sum_{t=1}^{T}\alpha_{t}%
V^{t}\right\Vert \leq\varepsilon\label{vApproximate}%
\end{equation}
for a small $\varepsilon>0,$ not true for the rest of the
eigenvalues.

We evaluate the approximation of each $\mathcal{L}$'s eigenvector $v$ using the 
eigenvectors of each time layer corresponding to the zero eigenvalue, $V^t$, by 
solving a regression problem where $\varepsilon_{i}$ is the $NT\times1$ vector of
residuals, and we denote the error at $i$ to be $\epsilon_{i}:=\left\Vert
\varepsilon_{i}\right\Vert$. Denote by $\lambda^{\ast}$ the first eigenvalue $\lambda_i$ for which $\epsilon_{i}>>\epsilon_{i-1}$.

\subsection{Discussion on the relation between edge density,
inter-layer weights and eigenvectors corresponding to the smallest
eigenvalues}

The present experimental results, in accordance with the developed theory, show that for a small eigenvalue of the supra-Laplacian $\mathcal{L}$, the eigenvectors $\psi^{R}$ and $\psi^{I}$ are approximations to the corresponding eigenvectors of the supra-Laplacian $\mathcal{L}.$ In Figure~\ref{fig:ER_eigen} we observe the eigenvectors of the supra-Laplacian of a temporal network composed of random Erdos-Renyi graphs,~\cite{Erdos1959}. The first few eigenvectors follow the same $\sin$ and $\cos$ functions as seen in Figure~\ref{fig:CosVsSin}, and thus can be used to identify the first order approximation by the constant block Jacobi model structure of the temporal network.

We investigate how the approximation of these eigenvectors is affected by the inter-layer weights and the density of the edge weights within each time layer. To showcase this, we simulate various benchmark temporal networks composed of random Erdos-Renyi networks with a varying degree of edge probability ${p=0.03,0.04,0.05,0.08,0.1,0.3}$ and inter-layer weights ${\omega=0.01,0.05,1,5}$, which are two factors that affect the approximation of the eigenvectors of the investigated supra-Laplacians $\mathcal{L}$, Figure~\ref{fig:err}.

Recall that we have denoted by $\lambda^{\ast}$  the smallest non-zero eigenvalue sensitive to within-layer connectivity patterns, i.e. breaking (\ref{vApproximate}). Then for all benchmark
networks types it is true that the value $\lambda^{\ast}$ is increasing
with a decreasing $\omega$ value: Smaller inter-layer weights $\omega$ lead to
greater separation between time layers, thus more eigenvectors behave as
predicted by perturbation theory. More eigenvectors are needed to explain each layer as separate. Higher inter-layer weights
influence more the resulting eigenvectors and fewer behave in a way as
predicted by perturbation theory. Lower inter-layer weights interfere less and
the behaviour of the eigenvectors resembles closely the behaviour of
eigenvectors as predicted by perturbation theory.

When the probability $p$ increases, the density within layers $A^{t}$ increases.
Since $\omega$ is fixed it cannot reflect on the increasing density of $A^{t}$
and the perturbation effect resulting from inter-layer matrices $W^{t,t+1}$ is
smaller. Thus for increasing $p$, i.e. for increasing density, the behaviour of more
eigenvectors resembles closely the behaviour of the eigenvectors as predicted
by perturbation theory.

When $p$ is decreasing, the eigenvalue $\lambda^{\ast}$ indicates that
more eigenvectors resemble closely the behaviour of eigenvectors as predicted
by perturbation theory. This is a result of the sparseness of the time layers
and the corresponding lower inter-layer weights $\omega_{i}^{t,t+1}$. The above observations
need further rigorous theoretical justification.

\begin{figure}
\centering
\includegraphics[width=.99\columnwidth]{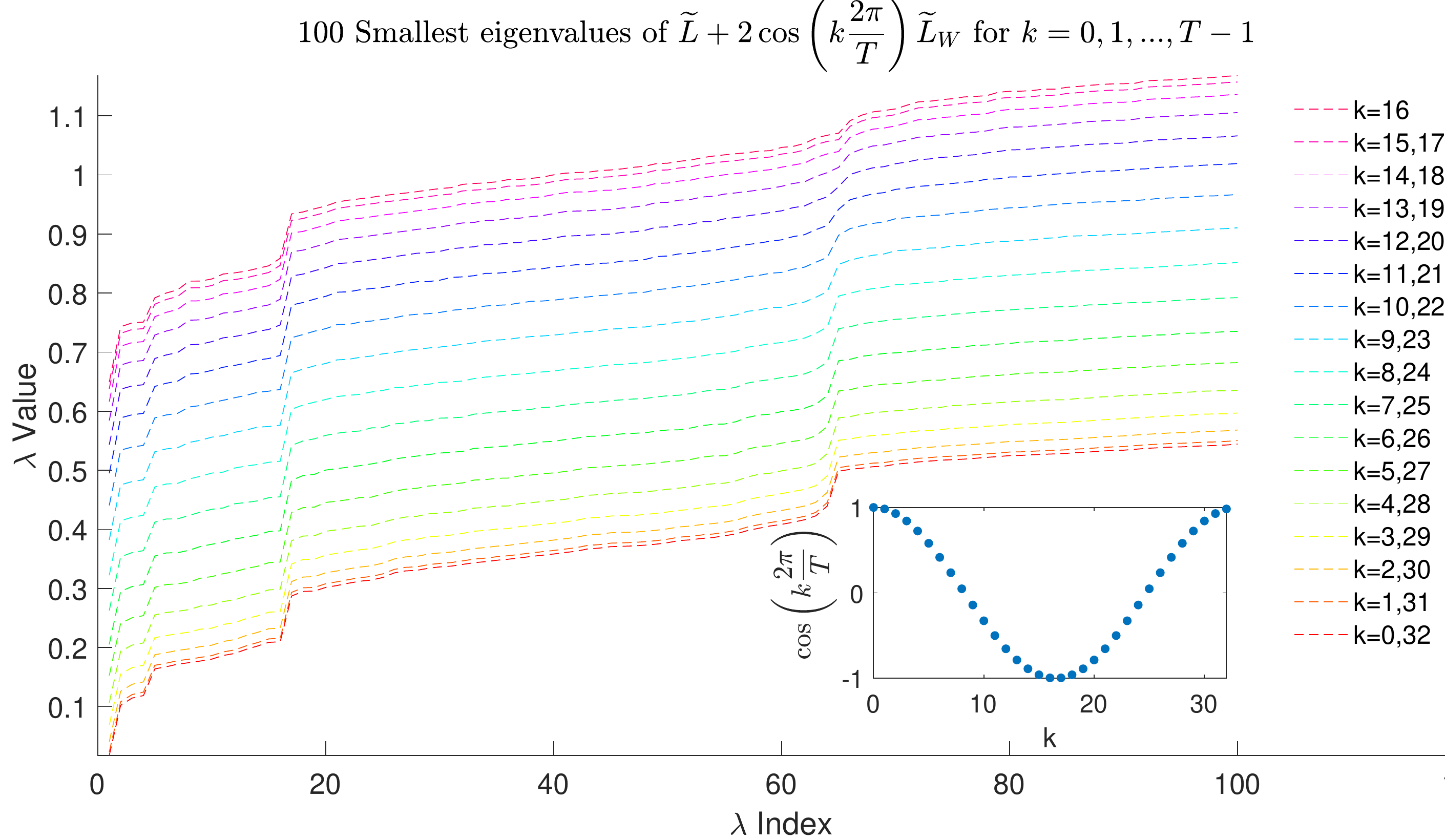}
\caption{\label{fig:eigenest_SP}\textbf{The \pmb{$100$} smallest eigenvalues of matrices \pmb{$\tilde{L}+2\cos\left(
k\frac{2\pi}{T}\right)  \tilde{L}_{W}$} for each \pmb{$k=0,1,2,...,32$}.} The matrices
$\tilde{L}$ and $\tilde{L}_{W}$ are obtained from a temporal network composed
of $T=33$ Sales-Pardo graphs each with $N=640$ nodes. The inter-layer weights $\omega$ are fixed at $1$. We include the
additional plot of $\cos\left(  k\frac{2\pi}{T}\right)  $ which determines the
monotonically increasing behaviour for eigenvalues for $0\leq k\leq15$ and
monotonically decreasing behaviour for eigenvalues for $17\leq k\leq32$.}%
\end{figure}

\subsection{Relation between the multi-scale community structure of the layers of a supra-Laplacian network and its eigenvalues.}
It is important to note that in Figure \ref{fig:eigenest_ER} the first few eigenvalues capture the block structure of the temporal network following the constant block Jacobi model, thus close to $0$, however after they start monotonically increasing without any clear cuts. From spectral graph partitioning~\cite{ding2001spectral} we know that this is indicative of the lack of structure within the networks, which is the case in here where each layer is a densely connected \textbf{Erdos-Renyi random graphs} with no community structure. In Figure~\ref{fig:eigenest_SP}, we demonstrate the behavior of the supra-Laplacian eigenvalues when each of the layers has multi-scale community structure simulated using the Sales-Pardo model,~\cite{Sales-Pardo2007}. Again the smallest eigenvalues capture the block structure of the temporal network, however, there are clear eigenvalue cuts where a new multi-scale community structure within the layers is captured. 

\section{Conclusions\label{sec6}}
The above results are crucial in interpreting spectral clustering properties of the supra-Laplacian matrix of all slowly-changing temporal networks that can be represented using a constant block Jacobi model. We have provided experimental results with Erdos-Renyi (unstructured) networks and Sales-Pardo hierarchical networks. Further investigation in these theoretical results will lead into more insights of the spectral properties of supra-Laplacian matrices for more general temporal networks. As presented in the paper, the above findings provide a fundamental understanding of the spectral properties of temporal networks on time periods where they are slowly changing which can significantly improve all spectral-based methods applied on temporal networks such as partitioning, node ranking, community detection, clustering, etc. The above results were successfully used to extend a multiscale community detection method,~\cite{Tremblay2014}, based on a spectral graph wavelets approach,~\cite{Hammond2009}, to temporal networks. The extended method,~\cite{Kuncheva2016}, takes advantage of the developed theory to automatically detect the different scales at which communities exist across layers, which is an advantage over the multilayer modularity maximization approach,~\cite{Mucha}, used for similar purposes. The above experimental results have been also replicated on temporal Sales-Pardo hierarchical benchmark networks, which are suitable for multi-scale community detection. There is also a detailed investigation of using inter-layer weights that account for the sparsity and similarity across layers,~\cite{KunchevaThesis}, including a real life application example to social networks data.
\section{Acknowledgements}
The author OK acknowledges the project KP-06-N52-1 with Bulgarian NSF. The author ZK acknowledges the project KP-06-N32-8 with Bulgarian NSF and EPSRC scholarship (2012-2016) at Imperial College London.
\bibliographystyle{unsrtnat}
\bibliography{LibNew}
\end{document}